\documentclass[10pt,twoside]{article}
\usepackage{graphicx}
\usepackage{amsmath}
\usepackage{Latex-document}

\newtheorem{definition}{Definition}[section]

\newtheorem{theorem}[definition]{Theorem}

\renewcommand{\thesection}{\arabic{section}}
\newcommand{\Section}[1]{\section{\hskip -1em.\hskip 0.7em#1}}

\def\res{\mathop{\hbox{\vrule height 7pt width .5pt
depth 0pt\vrule height .5pt width 6pt depth 0pt}}\nolimits}
\def\R{{\bf R}}

\title{\bf Optimal Transport Maps in\vskip -2mm Monge-Kantorovich Problem\vskip 6mm}

\author{L. Ambrosio\vspace*{-0.5cm}\thanks{Scuola Normale Superiore,
Piazza Cavalieri 7, 56126 Pisa, Italy. E-mail: luigi@ambrosio.sns.it}}
\date{\vspace{-8mm}}

\begin{document}
\maketitle

\thispagestyle{first} \setcounter{page}{131}

\begin{abstract}\vskip 3mm
In the first part of the paper we briefly decribe the classical
problem, raised by Monge in 1781, of optimal transportation of
mass. We discuss also Kantorovich's weak solution of the problem,
which leads to general existence results, to a dual formulation,
and to necessary and sufficient optimality conditions.

In the second part we describe some recent progress
on the problem of the existence of optimal transport maps. We show
that in several cases optimal transport maps can be obtained by
a singular perturbation technique based on the theory of
$\Gamma$-convergence, which yields as a byproduct
existence and stability results for classical
Monge solutions.

\vskip 4.5mm


\end{abstract}

\vskip 12mm

\Section{The optimal transport problem and its weak formulation}

\vskip-5mm \hspace{5mm}

In 1781, G.Monge raised in \cite{monge} the problem of
transporting a given distribution of matter (a pile of sand for instance)
into another (an excavation for instance) in such a way that
the work done is minimal. Denoting by $h_0,\,h_1:{\bf R}^2\to [0,+\infty)$
the Borel functions describing the initial and final distribution of
matter, there is obviously a compatibility condition, that the
total mass is the same:
\begin{equation}\label{balance}
\int_{\R^2}h_0(x)\,dx=\int_{\R^2}h_1(y)\,dy.
\end{equation}
Assuming with no loss of generality that the total mass is
$1$, we say that a Borel map $\psi:{\bf R}^2\to {\bf R}^2$ is a
{\em transport} if a local version of the balance of mass condition
holds, namely
\begin{equation}\label{balance1}
\int_{\psi^{-1}(E)}h_0(x)\,dx=\int_E h_1(y)\,dy
\qquad\mbox{\rm for any $E\subset\R^2$ Borel.}
\end{equation}
Then, the Monge problem consists in minimizing the work of
transportation in the class of transports, i.e.
\begin{equation}\label{Monge}
\inf\left\{\int_{\R^2}|\psi(x)-x|h_0(x)\,dx:\
\mbox{\rm $\psi$ transport}\right\}.
\end{equation}

The Monge transport problem can be easily generalized in many
directions, and all these generalizations have proved to be quite useful:

\noindent $\bullet$
General measurable spaces $X$, $Y$, with measurable maps $\psi:X\to Y$;

\noindent $\bullet$
General probability measures $\mu$ in $X$ and $\nu$ in $Y$.
In this case the local balance of mass condition (\ref{balance1})
reads as follows:
\begin{equation}\label{balance3}
\nu(E)=\mu(\psi^{-1}(E))\qquad\mbox{\rm for any $E\subset Y$ measurable.}
\end{equation}
This means that the push-forward operator $\psi_\#$ induced by $\psi$,
mapping probability measures in $X$ into probability measures in $Y$,
maps $\mu$ into $\nu$.

\noindent $\bullet$
General cost functions: a measurable map $c:X\times Y\to [0,+\infty]$.
In this case the cost to be minimized is
$$
W(\psi):=\int_X c\left(x,\psi(x)\right)\,d\mu(x).
$$

Even in Euclidean spaces, the problem of existence of optimal
transport maps is far from being trivial, mainly due to the
non-linearity with respect to $\psi$  of the condition $\psi_\#\mu=\nu$.
In particular the class of transports is not closed with
respect to any reasonable weak topology. Furthermore, it
is easy to build examples where the Monge problem is ill-posed
simply because there is no transport map: this happens for
instance when $\mu$ is a Dirac mass and $\nu$ is not a Dirac mass.

In order to overcome these difficulties, in 1942 L.V.Kantorovich
proposed in \cite{kant1} a notion of weak solution of the transport problem.
He suggested to look for {\em plans} instead of transports,
i.e. probability measures $\gamma$ in
$X\times Y$ whose marginals are $\mu$ and $\nu$. Formally this
means that $\pi_{X\#}\gamma=\mu$ and $\pi_{Y\#}\gamma=\nu$, where
$\pi_X:X\times Y\to X$ and $\pi_Y:X\times Y\to Y$ are the
canonical projections. Denoting by $\Pi(\mu,\nu)$ the
class of plans, he wrote the following minimization
problem
\begin{equation}\label{Kant}
\min\left\{\int_{X\times Y}c(x,y)\, d\gamma:\
\gamma\in\Pi(\mu,\nu)\right\}.
\end{equation}
Notice that $\Pi(\mu,\nu)$ is not empty, as the product $\mu\otimes\nu$
has $\mu$ and $\nu$ as marginals.
Due to the convexity of the new constraint $\gamma\in\Pi(\mu,\nu)$
it turns out that weak topologies can be effectively used to provide
existence of solutions to (\ref{Kant}): this happens for instance whenever
$X$ and $Y$ are Polish spaces and $c$ is lower semicontinuous
(see for instance \cite{RR}). Notice also that, by convexity of the energy,
the infimum is attained on a \emph{extremal} element of $\Pi(\mu,\nu)$.

The connection between the Kantorovich formulation of the transport
problem and Monge's original one can be seen noticing that any
transport map $\psi$ induces a planning $\gamma$, defined
by $(Id\times\psi)_\#\mu$. This planning is concentrated on the
graph of $\psi$ in $X\times Y$ and it is easy to show that the converse
holds, i.e. whenever $\gamma$ is concentrated on a graph, then
$\gamma$ is induced by a transport map. Since any transport induces
a planning with the same cost, it turns out that
$$
\inf\hbox{\rm (\ref{Monge})}\geq\min\hbox{\rm (\ref{Kant})}.
$$
Moreover, by approximating any plan by plans induced by
transports, it can be shown that equality holds under fairly
general assumptions (see for instance \cite{am1}).
Therefore we can really consider the Kantorovich formulation of
the transport problem as a weak formulation of the original problem.

If all extremal points of $\Pi(\mu,\nu)$ were
induced by transports one would get existence of
transport maps directly from the Kantorovich formulation.
It is not difficult to show that plannings $\gamma$ induced
by transports are extremal in $\Pi(\mu,\nu)$.
The converse holds in some very particular cases,
but unfortunately it is not true in general.
It turns out that the existence of optimal transport maps depends not
only on the geometry of $\Pi(\mu,\nu)$, but also (in a quite sensible way)
on the choice of the cost function $c$.

\Section{Existence of optimal transport maps}

\vskip-5mm \hspace{5mm}

In this section we focus on the problem of the existence of
optimal transport maps in the sense of Monge.
Before discussing in detail in the next sections
the two model cases in which the cost function is
the square of a distance or a distance (we refer to \cite{gangbocann}
for the case of \emph{concave} functions of the distance, not discussed
here), it is better to give an informal
description of the tools by now available for proving the existence
of optimal transport maps.

\smallskip
\noindent {\bf Strategy A} (Dual formulation). This strategy is based on the duality formula
\begin{equation}\label{duality}
\min \,(\hbox{\rm MK})=\sup\left\{\int_X h\,d\mu+\int_Yk\,d\nu
\right\},
\end{equation}
where the supremum runs among all pairs $(h,k)\in L^1(\mu)\times L^1(\nu)$
such that $h(x)+k(y)\leq c(x,y)$. The duality approach to the (MK) problem
was developed by Kantorovich, and then extended to more general cost
functions (see \cite{Kellerer}).
The transport map is obtained from an optimal pair $(h,k)$
in the dual formulation by making a first variation.
This strategy for proving the existence of an optimal
transport map goes back to the papers \cite{gangbo} and
\cite{caffarelli}.

\smallskip
\noindent {\bf Strategy B} (Cyclical monotonicity). In some situations the necessary (and sufficient) minimality
conditions for the primal problem, based upon the so-called $c$-cyclical monotonicity (\cite{smithknott1},
\cite{RR}, \cite{ruschendorf}) yield that \emph{any} optimal Kantorovich solution $\gamma$ is concentrated on a
graph $\Gamma$ (i.e. for $\mu$-a.e. $x$ there exists a \emph{unique} $y$ such that $(x,y)\in\Gamma$) and therefore
is
induced by a transport $\psi$.\\

This happens for instance when $c(x,y)=H(x-y)$, with $H$
\emph{strictly} convex in ${\bf R}^n$. This approach is pursued
in the papers \cite{gangbocann}, \cite{ruschendorf1}.

\smallskip
\noindent {\bf Strategy C} (Singular perturbation with strictly convex costs). One can try to get an optimal
transport map by making the cost strictly convex through a perturbation and then passing to the limit (see
\cite{caffelcann} and Theorem~\ref{tcaff}, Theorem~\ref{tcaff1} below). The main difficulty is to show (strong)
convergence at the level of the transport maps and not only at the level of transport plans.

\smallskip
\noindent {\bf Strategy D} (Reduction to a lower dimensional problem). This strategy has been initiated by
V.N.Sudakov in \cite{sudakov}. It consists in writing (typically through a disintegration) $\mu$ and $\nu$ as the
superposition of measures concentrated on lower dimensional sets and in solving the lower dimensional transport
problems, trying in the end to ``glue'' all the partial transport maps into a single transport map. This strategy
is discussed in detail in \cite{am1} and used, together with a ``variational'' decomposition, in \cite{am2}. The
simplest case is when the lower dimensional problems are $1$-dimensional, since the solution of the
$1$-dimensional transport problem is simply given by an increasing rearrangement, at least for convex functions of
the distance (see for instance \cite{am0}, \cite{RR}, \cite{villani}).

\smallskip
Strategies A and B are basically equivalent and yield existence and
uniqueness at the same time: the first one could be
preferable for someone, as a very small measure-theoretic apparatus
is involved. On the other hand, it strongly depends on the existence
of maximizing pairs in the dual formulation, and this existence issue
can be more subtle than the existence issue for the primal problem
(see \cite{RR} and the discussion in \cite{am1}).
For this reason it seems that the second strategy
can work for more general classes of cost functions.

Strategies C and D have been devised to deal with situations where the cost
function is convex but not strictly convex. Also these two strategies are
closely related, as the strictly convex perturbation often leads to an
effective dimension reduction of the problem (see for instance \cite{am2}).

\Section{cost=distance$^2$}

\vskip-5mm \hspace{5mm}

In this section we consider the case when $X=Y$ and the cost function $c$
is proportional to the square of a distance $d$. For convenience we
normalize $c$ so that $c=d^2/2$. The first result in the Euclidean space
${\bf R}^n$ has been discovered independently by many authors
Y.Brenier \cite{brenier}, \cite{brenier1},
S.T.Rachev and L.R.\"uschendorf \cite{rr1}, \cite{ruschendorf},
and C.Smith and M.Knott \cite{smithknott}.

\begin{theorem}\label{tbrenier}
Assume that $\mu$ is absolutely continuous with respect to
${\cal L}^n$ and that $\mu$ and $\nu$ have finite second order moments. Then
there exists a unique optimal transport map $\psi$. Moreover
$\psi$ is the gradient of a convex function.
\end{theorem}

In this case the proof comes from the fact that both strategies A and B yield that
the displacement $x-\psi(x)$ is the gradient of a $c$-concave function, i.e. a function
representable as
$$
h(x)=\inf_{(y,t)\in I} c(x,y)+t\qquad\forall x\in {\bf R}^n
$$
for a suitable non-empty set $I\subset Y\times{\bf R}$.
The concept of $c$-concavity \cite{ruschendorf} has been extensively used
to develop a very general duality theory for the (MK) problem, based on (\ref{duality}). In this
special Euclidean situation it is immediate to realize that $c$-concavity of $h$ is equivalent
to concavity (in the classical sense) of $h-\tfrac{1}{2}|x|^2$, hence
$$
\psi(x)=x-\nabla h(x)=\nabla \left[\frac{1}{2}|x|^2-h(x)\right]
$$
is the gradient of a convex function.
Finally, notice that the assumption on $\mu$ can be sharpened
(see \cite{gangbocann}), assuming for
instance that $\mu(B)=0$ whenever $B$ has finite
${\cal H}^{n-1}$-measure.
This is due to the fact that the non-differentiability set of a concave function
is $\sigma$-finite with respect to ${\cal H}^{n-1}$ (see for instance
\cite{alam}). Also the assumption
about second order moments can be relaxed, assuming only that the infimum
of the (MK) problem with data $\mu,\,\nu$ is finite.

The following result, due to R.Mc Cann \cite{mccann}, is much more recent.

\begin{theorem}
Assume that $M$ is a $C^3$, complete Riemannian manifold with no
boundary and $d$ is the Riemannian distance. If $\mu$, $\nu$ have
finite second order moments and $\mu$ is absolutely continuous
with respect to ${\rm vol}_M$ there exists a unique optimal
transport map $\psi$.

Moreover there exists a $c$-concave potential $h:M\to{\bf R}$
such that
$$
\psi(x)={\rm exp}_x\left(-\nabla h(x)\right)\quad \hbox{\rm ${\rm vol}_M$-a.e..}
$$
\end{theorem}

This Riemannian extension of Theorem~\ref{tbrenier} is non trivial, due to the fact that $d^2$ is not smooth in
the large. The proof uses some semiconcavity estimates for $d^2$ and the fact that $d^2$ is $C^2$ for $x$ close to
$y$ (this is where the $C^3$ assumption on $M$ is needed). It is interesting to notice that the results of
\cite{manmen} (where the eikonal equation is read in local coordinates), based on the theory of viscosity
solutions --- see in particular Theorem~5.3 of \cite{pllions} --- allow to push Mc Cann's technique up to $C^2$
manifolds.

Can we go beyond Riemannian manifolds in the existence theory? A model case is given by stratified Carnot groups
endowed with the Carnot-Carath\'eodory metric $d_{CC}$, as these spaces arise in a very natural way as limits of
Riemannian manifolds with respect to the Gromov-Hausdorff convergence (see \cite{gromov}). At this moment a
general strategy is still missing, but some preliminary investigations in the Heisenberg group $H_n$ show that
positive results analogous to the Riemannian ones can be expected. The following result is proved in \cite{am3}:

\begin{theorem}
If $n=1,\,2$ and $\mu$ is a probability measure in $H_n$
absolutely continuous with respect to ${\cal L}^{2n+1}$, then:\\
(a) there exists a unique optimal transport map $\psi$, deriving from
a $c$-concave potential $h$;\\
(b) If $d_p\uparrow d_{CC}$ are Riemannian left invariant metrics
then Mc Cann's optimal  transport maps $\psi_p$ relative to
$c_p=d_p^2/2$ converge in measure to $\psi$ as $p\to\infty$.
\end{theorem}

The restriction to $H_n$, $n\leq 2$, arises from
the fact that so far we have been able to carry on some explicit
computations only for $n\leq 2$. We expect that this restriction
could be removed. The proof of (b) is not direct,
as Mc Cann's exponential representation
$\psi_p={\rm exp}^p_x(-\nabla^p h_p)$
``degenerates" as $p\to\infty$, because the injectivity
radius of the approximating manifolds tends to $0$. This is due to
the fact that in CC metric spaces geodesics
exist but are not unique, not even in the small.

Finally, if we replace $c$ by the square of the Kor\'anyi norm (related to the fundamental solution of the Kohn
sub-Laplacian), namely
$$
\tilde{c}(x,y):=\frac{1}{2}\Vert y^{-1}x\Vert^2\quad\hbox{\rm with}\quad
\Vert (z,t)\Vert:={^{4}}\sqrt{|z|^4+t^2}
$$
(here we identify $H_n$ with ${\bf C}^n\times\R$)
then we are still able to prove existence in any Heisenberg group
$H_n$. The proof uses some fine properties of $BV$ functions on
sub-Riemannian groups \cite{amagnani}. However, we can't hope for a
Riemannian approximation
result, as the Kor\'anyi norm induces a metric $d_K$ which is not
geodesic. It turns out that the geodesic metric associated to $d_K$
is a constant multiple of $d_{CC}$.

\Section{cost=distance}

\vskip-5mm \hspace{5mm}

In this section we consider the case when $X=Y$ and the cost function $c$
is a distance. In this case both strategies A and B give only a partial
information about the location of $y$, for given $x$. In particular it
is not true that any optimal Kantorovich plan $\gamma$ is induced by a transport map.
Indeed, if the first order moments of $\mu$ and $\nu$ are finite,
the dual formulation provides us with a maximizing pair $(h,k)=(u,-u)$,
with $u:X\to{\bf R}$ $1$-Lipschitz. If $X={\bf R}^n$ and the distance is
induced by a norm $\Vert\cdot\Vert$, this provides the implication
\begin{equation}\label{nonno}
(x,y)\in{\rm spt\,}\gamma\qquad\Longrightarrow\qquad
y\in\left\{x-s\xi:\ \xi\in\left(du(x)\right)^*,\,\,\,s\geq 0\right\}
\end{equation}
at any differentiability point of $u$. Here we consider the natural duality map
between covectors and vectors given by
$$
L^*:=\left\{\xi\in{\bf R}^n:\ L(\xi)=\Vert L\Vert_*\,\,\,\hbox{\rm and}\,\,\,\Vert\xi\Vert=1\right\}.
$$
The most favourable case is when the norm is strictly convex (e.g. the Euclidean norm):
in this situation the $*$ operator is single-valued and we recover from (\ref{nonno}) an information
on the direction of transportation, i.e. $\left(du(x)\right)^*$,
but not on the length of transportation. If the norm is
not strictly convex (e.g. the $l_1$ or $l_\infty$ norm) then even the information on the
direction of transportation, encoded in $\left(du(x)\right)^*$, is partial.

The first attempt to bypass these difficulties came with the work
of V.N.Sudakov \cite{sudakov}, who claimed to have a solution for any
distance cost function induced by a norm. Sudakov's approach is based on
a clever decomposition of the space $\R^n$ in affine regions with variable
dimension where the Kantorovich dual potential $u$ associated to the transport
problem is an affine function. His strategy is to solve the
transport problem in any of these regions, eventually getting
an optimal transport map just by gluing all these transport
maps. An essential ingredient in his proof is Proposition~78,
where he states that, if $\mu<<{\cal L}^n$, then the conditional
measures induced by the decomposition are absolutely continuous
with respect to the Lebesgue measure (of the correct dimension).
However, it turns out that this property is not true in general
even for the simplest decomposition, i.e. the decomposition in segments:
G.Alberti, B.Kirchheim and D.Preiss found an example
of a compact
faily of pairwise disjoint open segments in $\R^3$ such that the family
$M$ of their midpoints has strictly positive Lebesgue measure
(the construction is a variant of previous examples due
to A.S.Besicovitch and D.G.Larman, see also \cite{am0} and \cite{am2}).
In this case, choosing $\mu={\cal L}^3\res M$, the conditional measures
induced by the decomposition are Dirac masses.
Therefore it is clear that this kind of counterexamples should
be ruled out by some kind of additional ``regularity'' property of
the decomposition. In this way the Sudakov strategy would be
fully rigorous. As noticed in \cite{am2}, this regularity comes
for free only in the case $n=2$, using the fact that transport rays
do not cross in their interior.

Several years later, L.C.Evans and W.Gangbo made a remarkable progress
in \cite{evansgangbo}, showing
by differential methods the existence of a transport map,
under the assumption that
${\rm spt\,}\mu\cap{\rm spt\,}\nu=\emptyset$, that the two
measures are absolutely continuous with respect to ${\cal L}^n$
and that their densities are Lipschitz functions with compact
support. The missing
piece of information about the length of transportation is
recovered by a $p$-laplacian approximation
$$
-{\rm div\,}\left(|\nabla u|^{p-2}\nabla u\right)=\mu-\nu,\qquad
u\in H^1_0(B_R),\qquad R\gg 1
$$
obtaining in the limit as $p\to +\infty$ a nonnegative
function $a\in L^\infty(\R^n)$
and a $1$-Lipschitz function $u$ solving
$$
-{\rm div\,}(a\nabla u)=\mu-\nu,\qquad
\hbox{\rm $|\nabla u|=1$ ${\cal L}^n$-a.e. on $\{a>0\}$.}
$$

The diffusion coefficient $a$ in the PDE above plays a special role in the
theory. Indeed, one can show (see \cite{am0}) that
the measure $\sigma:=a{\cal L}^n$, the so-called {\em transport density},
can be represented in several different way, and in particular as
\begin{equation}\label{freeze}
\sigma(B)=\int {\cal H}^1\left(B\cap [x,y]\right)\,d\gamma(x,y)
\qquad\forall\hbox{\rm $B\subset {\bf R}^n$ Borel}
\end{equation}
for some optimal planning $\gamma$. Notice that the total mass of
$\sigma$ is $\int|x-y|\,d\gamma$, the total work done and the meaning
of $\sigma(B)$ is the work done \emph{within} $B$ during the transport
process. This representation of the transport density has been introduced
by G.Bouchitt\'e and G.Buttazzo in \cite{buttazzo}, who showed that the
a constant multiple of the transport density is a solution of their
so-called mass optimization
problem. Later, in \cite{am0}, it was shown that there is actually
a 1-1 correspondence between solutions of the mass optimization
problem and transport densities, defined as in (\ref{freeze}).

One can also show (\cite{am0}, \cite{depapra}, \cite{felcann},
\cite{depapraevans}) that $\sigma$ is unique (unlike $\gamma$)
if either $\mu$ or $\nu$ are absolutely continuous.
Moreover, the nonlinear operator mapping
$(\mu,\nu)\in L^1\times L^1$ into $a\in L^1$ maps $L^p\times L^p$
into $L^p$ for $1\leq p\leq\infty$.

Coming back to the problem of the existence of optimal
transport maps with Euclidean distance $|x-y|$
(or, more generally, with a distance induced by a $C^2$ and
uniformly convex norm), the first existence results for
general absolutely continuous measures $\mu,\,\nu$
with compact support have been
independently obtained by L.Caffarelli, M.Feldman and R.Mc Cann in
\cite{caffelcann} and by N.Trudinger and L.Wang in \cite{trudinger}.
Afterwards, the author estabilished in \cite{am0} the existence
of an optimal transport map assuming only that the initial measure
$\mu$ is absolutely continuous, and the results of \cite{caffelcann}
and \cite{trudinger} have been extended to a Riemannian setting
in \cite{feldman1}.
All these proofs involve basically
a Sudakov decomposition in transport rays, but the technical
implementation of the idea is different from paper to paper:
for instance in \cite{caffelcann}
a local change of variable is made, so that transport rays become
parallel and Fubini theorem, in place of abstract disintegration
theorems for measures, can be used. The proof in \cite{am1}, instead,
uses the co-area formula to show that absolute continuity with respect
to Lebesgue measure is stable under disintegration.

The following result \cite{am1} is a slight improvement of \cite{caffelcann},
where existence of an optimal transport map was estabilished but not
the stability property. The result holds under regularity and uniform
convexity assumptions for the norm $\Vert\cdot\Vert$.

\begin{theorem}\label{tcaff}
Let $\mu,\nu$ be with compact support, with $\mu<<{\cal L}^n$,
and let $\psi_\epsilon$ be the unique optimal transport maps
relative to the costs $c_\epsilon(x,y):=\Vert x-y\Vert^{1+\epsilon}$.
Then $\psi_\epsilon$ converge as $\epsilon\downarrow 0$ to an
optimal transport map $\psi$ for $c(x,y)=\Vert x-y\Vert$.
\end{theorem}

The proof is based only the fact that any plan $\gamma_0$, limit of
some sequence of plans $(Id\times\psi_{\epsilon_i})$, is not only
optimal for the (MK) problem, but also for the \emph{secondary} one
\begin{equation}\label{corea}
\min_{\gamma\in\Pi_1(\mu,\nu)}
\int_{\R^n\times\R^n}\Vert x-y\Vert\ln(\Vert x-y\Vert)\,d\gamma,
\end{equation}
where $\Pi_1(\mu,\nu)$ denotes the class of all optimal
plannings for the Kantorovich problem
(the entropy function in (\ref{corea}) comes from
the Taylor expansion of $c_\epsilon$ around $\epsilon=0$).
It turns out that this additional minimality property selects a
unique plan induced by a transport $\psi$ and, a posteriori,
$\psi$ is the same
map built in \cite{caffelcann}. A class of counterexamples built in
\cite{am1} shows that the absolute continuity assumption on $\mu$
cannot be weakened, unlike the strictly convex case.

This ``variational'' procedure seems to select extremal elements of
$\Pi(\mu,\nu)$ in a very effective way. This phenomenon is apparent
in view of the following result \cite{am2}, which holds for all
``crystalline'' norms $\Vert\cdot\Vert$
(i.e. norms whose unit sphere is contained
in finitely many hyperplanes).

\begin{theorem}\label{tcaff1}
Let $\mu,\nu$ be as in Theorem~\ref{tcaff} and let
$\psi_\epsilon$ be the unique optimal transport maps
relative to the costs
$$
c_\epsilon(x,y):=\Vert x-y\Vert+\epsilon|x-y|+
\epsilon^2|x-y|\ln |x-y|.
$$
Then $\psi_\epsilon$ converge as $\epsilon\downarrow 0$ to an
optimal transport map $\psi$ for $c(x,y)=\Vert x-y\Vert$.
\end{theorem}

In this case a secondary and a ternary variational problem are
involved, and we show that the latter has a unique solution which is
also induced by a transport.

Some borderline cases between ``crystalline'' norms and
``Euclidean'' norms apparently can't be attacked by any of the
existing techniques. In particular the existence of optimal
transport maps for the cost induced by a general norm in ${\bf R}^n$,
$n\geq 3$, is still open.

\label{lastpage}

\end{document}